\documentclass[12pt]{amsart}   
\usepackage{amssymb}
\usepackage{amsmath}
\usepackage{graphicx}      

\usepackage{amscd}
\newtheorem{theorem}{Theorem}
\theoremstyle{plain}

\newtheorem{result}[theorem]{Result}

\newtheorem{claim}[theorem]{Claim}

\newtheorem{conjecture}[theorem]{Conjecture}

\newtheorem{definition}[theorem]{Definition}

\newtheorem{lemma}[theorem]{Lemma}
\newtheorem{notation}[theorem]{Notation}

\input{tcilatex}

\begin{document}
\title{The pebbling threshold of the square of cliques}
\author{A. Bekmetjev}
\address{A. Bekmetjev\\ Department of Mathematics \\ Hope College\\ Holland, MI 49422-9000}
\email{bekmetjev@hope.edu}
\author{G. Hurlbert \thanks{}}
\address{G. Hurlbert\\ Department of Mathematics and Statistics \\ Arizona State University\\ Tempe, AZ 85287-1804}
\email{hurlbert@asu.edu}

\thanks{$^*$ Partially supported by National Security Agency grant 
\#MDA9040210095.}
\date{December 30, 2003}
\maketitle
\begin{abstract}
Given an initial configuration of pebbles on a graph, one can move pebbles
in pairs along edges, at the cost of one of the pebbles moved, with the 
objective of reaching a specified target vertex.
The pebbling number of a graph is the minimum number of pebbles so that
every configuration of that many pebbles can reach any chosen target.
The pebbling threshold of a sequence of graphs is roughly the number of
pebbles so that almost every (resp. almost no) configuration of asymptotically 
more (resp. fewer) pebbles can reach any chosen target.
In this paper we find the pebbling threshold of the sequence of
squares of cliques, improving upon an earlier result of Boyle and
verifying an important instance of a probabilistic version of
Graham's product conjecture.
\end{abstract}



\section{Pebbling Number}

Consider a connected graph $G$ on $n$ vertices. Suppose that a 
\emph{configuration} $C$ of $t$ pebbles is
placed onto the vertices of graph $G.$ 
A\emph{\ pebbling step\ }from $u$ to $v$ consists of removing two
pebbles from vertex $u$\ and then placing one pebble on an adjacent vertex $%
v $\emph{. }We say that a pebble can be \textit{moved} to a vertex $r$
(called \textit{root} vertex) if after finitely many steps $r$ has at least
one pebble.\emph{\ }If it is possible to move a pebble to the root vertex $r$%
\ then we say that $C$\ is $r$\emph{-solvable}; otherwise, $C$\ is $r$\emph{%
-unsolvable}. Finally, we call $C$\ \emph{solvable}\ if it is $r$-solvable
for all $r$, and \emph{unsolvable}\ otherwise. Define
the \emph{pebbling number}\ $\pi (G)$\ to be the smallest integer $t$\ such
that every configuration of $t$\ pebbles on the vertices of $G$\ is solvable.
A fair amount is known about the pebbling numbers of typical graphs like
complete graphs, paths, cycles, cubes, etc. (see \cite{hurl} for a survey),
relations to known parameters such as connectivity \cite{chkt}, 
diameter \cite{bukh}, girth \cite{ch}, and domination number \cite{changod},
and interesting variations such as optimal pebbling \cite{moewsopt}, and
cover pebbling \cite{ccfhpst} are being investigated.

\section{Random Configurations}

In this paper we consider a random pebbling model in which a particular
configuration of pebbles is selected uniformly at random from the set of all
configurations with a fixed number of pebbles. One can think of the configuration of pebbles as a placement of unlabeled balls in
labeled distinct urns. This is analogous to the so-called \emph{static} model
of random graphs, whose sample space consists of all graphs with a fixed
number of edges. Since vertices may have more than one
pebble, a particular configuration is a multiset of $t$ elements with the
ground set $[n]$. We construct the probability space $\mathbf{C}_{n,t}$ by
choosing configurations randomly and assuming that they all are equally
likely to occur.

The size of the set $\mathbf{C}_{n,t}$ is the number of possible
arrangements of $t$ identical balls placed in $n$ distinct urns, so $|%
\mathbf{C}_{n,t}|={\binom{n+t-1}{t}}$, which we denote by
$\QATOPD\langle \rangle {n}{t}$ (the reader may find it useful to use the
terminology ``$n$ \textit{pebble} $t$''). We will be interested in the
probability spaces associated with sequences of graphs $\mathcal{G}%
=(G_{1},G_{2},\ldots ,G_{n},\ldots ).$ In this notation the index $n$
represents the position of the graph $G_{n}.$ In some of the graph
sequences, such as cubes, for example, the size of the vertex set of $G_{n}$
is not the same as the position. Therefore, we define $N=N_n=N(G_n)=|V(G)|$ 
to be the number of vertices of $G_n$. Graphs in $\mathcal{G}$ are in
ascending order with respect to this number, i.e. $N_n>N_m$ for $n>m$.

We will study the pebbling threshold phenomenon that occurs in this model,
as it does for many random graph properties. For two functions $f(n)$\ and $%
g(n)$\ we write $f\ll g$, (equivalently $g\gg f$\ ) if 
$\lim_{n\rightarrow\infty} f(n)/g(n)=0$.
We set $o(g)=\{f$\ $|$\ $f\ll g\}$\ and $%
\omega (f)=\{g$\ $|$\ $f\ll g\}$. Also, we write $f\in O(g)$, or
equivalently $g\in \Omega (f),$ when there are positive constants $c$\ and $%
k$\ such that $f(n)/g(n)<c$, for all $n>k$.\ In particular, if$\
f(n)/g(n)\rightarrow 1$ as $n\rightarrow \infty$, we write $%
f\thicksim g$. Furthermore, we define $\Theta (g)$\ $=$\ $O(g)\cap \Omega
(g)$. Finally, for two sets of functions $F$ and $G$\ we write $F\lesssim G$%
\ if $f\in O(g)$ for all $f\in F,g\in G$\textbf{.}

A function $f=f(n)$\ is called a\emph{\ threshold\ }for the graph sequence $%
\mathcal{G}$, and we write\emph{\ }$f\in \tau _{\mathcal{G}}$\emph{, }if 
$P_{\mathcal{G}}(n,t)\rightarrow 1$ whenever $t\gg f,$ and $P_{\mathcal{G}%
}(n,t)\rightarrow 0,$ whenever\emph{\ }$t\ll f.$\emph{\ }
In other words, if $f=f(n)\in \tau _{\mathcal{G}}$,
then for any function $\varpi
=\varpi (n)$ tending to infinity with $n$, $P_{\mathcal{G}}(n,\varpi
f)\rightarrow 1$ and $P_{\mathcal{G}}(n,f/\varpi )\rightarrow 0$ as\emph{\ }%
$n\rightarrow \infty $.  

Roughly speaking, the pebbling number describes the ``worst-case'' scenario,
as it is one more than the size of the largest unsolvable configuration.
The threshold function, on the
other hand, deals with ``typical'' configurations and estimates the
average chance of being solvable. For example, the
threshold of family of cliques $\mathcal{K}$ is $\tau _{\mathcal{K}%
}=\Theta (\sqrt{N}).$ This problem is similar to the well-known ``birthday'' 
problem --- 
how many people must be in a room so that with high probability two people
share the same birth date?
--- but here the pebbles are unlabelled.  The
general existence of the pebbling threshold is established in \cite{thres},
and in \cite{cehk} it is shown that every graph sequence $\mathcal G$ satisfies
$\tau_{\mathcal G}\subset \Omega(f)\cap O(g)$, where $f\in\tau_{\mathcal K}$ 
and $g\in\tau_{\mathcal P}$, for the sequence of paths ${\mathcal P}$.
We are going to compute the pebbling threshold of the sequence of
squares of cliques, thereby verifying an instance of the threshold analogue 
of Graham's product conjecture.

\section{Cartesian Products and Graham's Conjecture}

Chung's paper \cite{chung} raised a natural question about the relationship
between the pebbling numbers of individual graphs and the pebbling number of
their cartesian product.

\begin{definition}
\emph{The} Cartesian product\emph{\ of two graphs }$G_{1}$\emph{\ and }$%
G_{2},$\emph{\ denoted }$G_{1}\Box G_{2}$\emph{\ is the graph with }$\emph{%
vertex\ set}$\emph{\ }
\begin{equation*}
V(G_{1}\Box G_{2})=\{(v_{1},v_{2})\text{ }|\text{ }v_{1}\in
V(G_{1}),v_{2}\in V(G_{2})\}
\end{equation*}
\emph{and }$\emph{edge\ set}$%
\begin{equation*}
\begin{tabular}{lll}
$E(G_{1}\Box G_{2})$ & $=$ & $\{$ $((v_{1},v_{2}),(w_{1},w_{2}))$ $|$ $%
v_{1}=w_{1}\text{\emph{ and} }(v_{2},w_{2})\in E(G_{2})$ \emph{or} \\ 
&  & $\qquad \qquad \quad \quad \qquad $\emph{\ \ \ \ \ }$v_{2}=w_{2}\text{%
\emph{ and }}(v_{1},w_{1})\in E(G_{1})$ $\}$\ .
\end{tabular}
\end{equation*}
\emph{\ }
\end{definition}

The general conjecture about the pebbling number of the cartesian product of
graphs was originally stated by Graham (\cite{chung}).
\begin{conjecture}
\label{graham} 
For all graphs $G_{1}$ and $G_{2}$ we have that 
\begin{equation*}
\pi (G_{1}\Box G_{2})\leqslant \pi (G_{1})\pi (G_{2})\ .
\end{equation*}
\end{conjecture}
There are several results supporting this conjecture. It is known \cite{chung} 
that the $m$-dimensional cube and that the product of cliques satisfy 
this conjecture.
Also, Moews \cite{moews} proved it holds for the product of trees.
Pachter et al. \cite{pachter} proved the conjecture for the product of
cycles with some exceptions: it holds for $C_{m}\Box C_{n}$ where $m$ and $n$
are not both from the set $\{5,7,9,11,13\}$. Herscovici and Higgins in \cite
{higgins} proved it for $C_{5}\Box C_{5}$. Recently, Herscovici \cite
{higginsnew} found a proof for all these exceptions confirming Graham's
conjecture for the product of cycles.
Finally, the conjecture holds for dense graphs \cite{ch}.

\section{\label{secprobgraham} Threshold Version and Main Theorem}

For the graph sequences ${\mathcal G}=(G_{1},\ldots G_{n},\ldots )$ and 
${\mathcal H}=(H_{1},\ldots ,H_{n},\ldots )$ let us define the sequence 
${\mathcal G}\Box{\mathcal H}=
(G_{1}\Box H_{1},\ldots, G_{n}\Box H_{n},\ldots )$. 
The sequence 
$\mathcal{G\Box H}$\ is called the \emph{cartesian product }of\emph{\ }
$\mathcal{G}$\emph{\ }and\emph{\ }$\mathcal{H}$. The number of
vertices of the $n^{\emph{th}}$ element of $\mathcal{G\Box H}$ is
$N(G_{n}\Box H_{n})=N(G_{n})N(H_{n})$.
Here we are interested in the following probabilistic version of
Conjecture \ref{graham}. 
\begin{conjecture}\label{hurlcon}
Let $\mathcal{F}$ and $\mathcal{G}$ be two graph
sequences with numbers of vertices $R=N(F_n)$ and $S=N(G_n)$, respectively,
and with pebbling thresholds $\tau_{\mathcal F}$ and $\tau_{\mathcal G}$,
respectively.
Let $f\in\tau_{\mathcal F}$, $g\in\tau_{\mathcal G}$, and 
$h\in\tau_{\mathcal H}$, where ${\mathcal H}={\mathcal F}\Box{\mathcal G}$,
having $T=N(H_n)=RS$ vertices.
Then
\begin{equation*}
h(T)\in O(f(R)g(S))\ .
\end{equation*}
\end{conjecture}

This conjecture is shown to hold for $d$-dimensional grids (products of
paths) in \cite{ch}.
We are going to verify Conjecture \ref{hurlcon} for the cartesian 
product of
cliques $\mathcal{K}^{2}=\mathcal{K\Box K}=(K_{1}\Box K_{1},\ldots,
K_{n}\Box K_{n},\ldots ).$ 
If true, the pebbling threshold for the product of cliques should be
\begin{equation*}
\tau _{\mathcal{K}^{2}}\subseteq \Theta (\sqrt{N^{1/2}}\sqrt{N^{1/2}}%
)=\Theta (\sqrt{N})\ ,
\end{equation*}
where $N$ is the number of vertices of $\mathcal{K}^{2},$ namely $N=n^{2}$.
This would improve Boyle's \cite{boyle} result that $\tau _{\mathcal{K}%
^{2}}\subseteq O(N^{3/4})$ and give the exact result (recall the
lower bound for all sequences mentioned above).
Our main result is the following theorem.
\begin{theorem}
\label{thmk2}Let $\mathcal{K}^{2}$ be the sequence of the cartesian products
of cliques, with $N=N(K_{n}^{2})$. Then the pebbling threshold of
$\mathcal{K}^{2}$ is 
\begin{equation*}
\tau _{\mathcal{K}^{2}}=\Theta (\sqrt{N})\ .
\end{equation*}
\end{theorem}
This theorem is perhaps surprising, considering that the graph $K_n^2$ is
fairly sparse.  It seems that the structure of the graph is what keeps its
threshold small.

\section{Cops and Robbers}

\begin{figure}[tbp]
\centering
\resizebox {8cm}{7cm}
{\includegraphics[0,0][281,249]{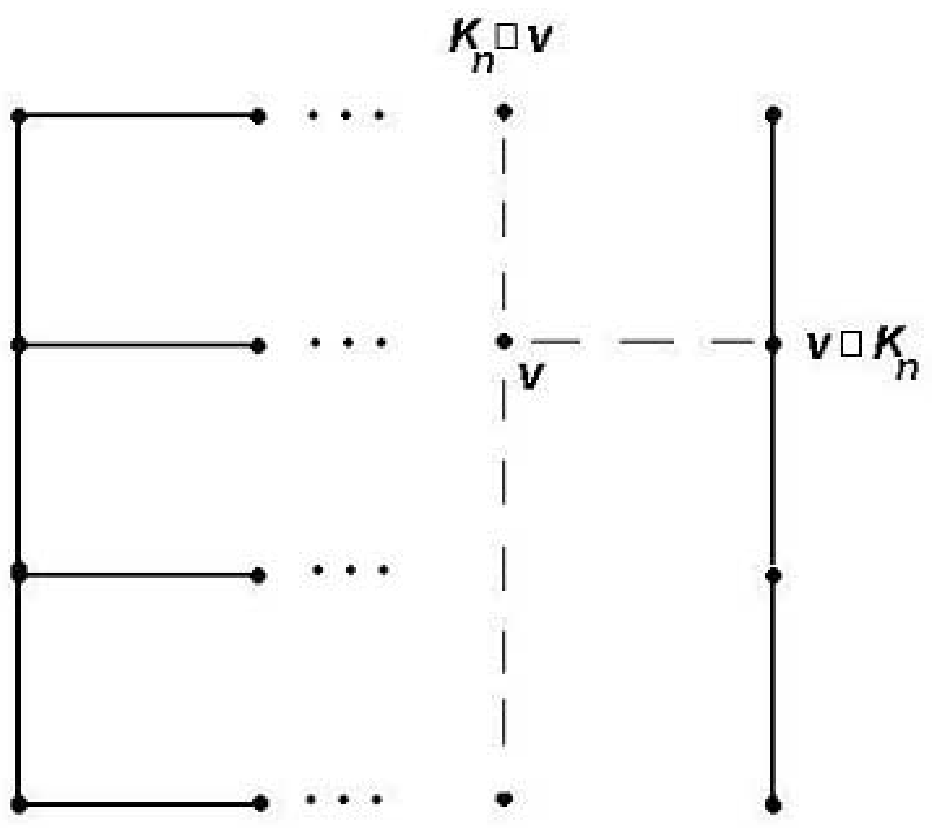} }
\caption{Schematic Presentation of $K_{n}\Box K_{n}$.}
\label{cliquebyclique}
\end{figure}

Let us consider a particular configuration of pebbles on the cartesian
product of cliques $K_{n}\Box K_{n}$. Note that this graph can be thought
of as a rectangular grid with each row and column a complete graph.
Therefore, to pebble to a specific root $r$ one needs to collect
two or more pebbles on any vertex that belongs to the \textit{row}
$r\Box K_{n}$ or to the \textit{column} $%
K_{n}\Box r$ (see Figure
\ref{cliquebyclique}). This suggests the following interpretation of the
pebbling problem.  We partition the vertices of $%
K_{n}^{2}$ into three
distinct sets: \textit{police}, or \textit{cops} ($\mathbf{P}$%
), \textit{citizens}
($\mathbf{T}$) and \textit{robbers} ($\mathbf{R}$).  Vertices in
the set $\mathbf{P}$ are those with two or more pebbles on them, $%
\mathbf{T}$ is the set of vertices with one pebble, and $\mathbf{R}$
is the remaining set of empty vertices.  (This approach is motivated by
a variety ``Cops and Robbers'' games, one of the more prevalent types of
games on graphs. More information on these types of games can be found,
for example, in \cite{adler, coprob2, copsandrobbers, seymour}.)

\begin{figure}[tbp]
\centering
\resizebox {8cm}{7cm}
{\includegraphics[0,0][330,311]{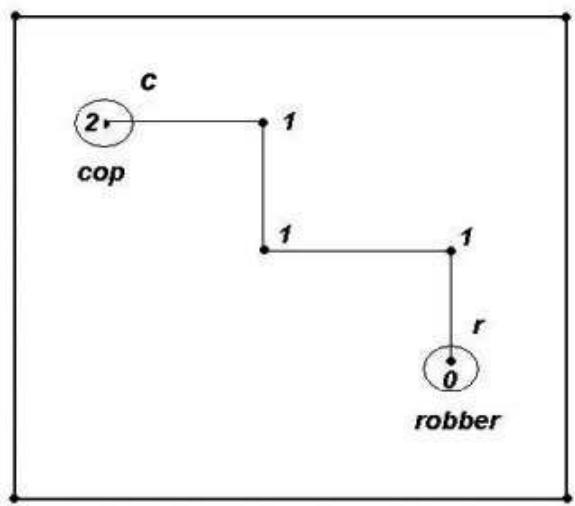} }
\caption{Pebbling from ``Cop'' to ``Robber''.}
\label{policepath}
\end{figure}

In our case, the robber is immobile and cops can move only in certain
directions and their number may change during the game. If root $r$ is
chosen in $\mathbf{R}$ then for a pebbling configuration to be solvable it
is sufficient that there is at least one cop 
on $r\Box K_{n}$ or $K_{n}\Box r$. Any citizen can become a cop
if it is possible to move at least one pebble to it from some other cop.
We say that vertex $u$ \textit{sees }$v$ if $u$ and $v$ are in the
same row or column of $K_{n}^{2}.$ Furthermore, we say that a robber $%
r=v_{0} $ can be \textit{caught }if there is a sequence of citizens $%
v_{1},\ldots ,v_{k-1}$ and a cop $c=v_{k}$ so that $v_{i}$ sees $v_{i+1}$
for $0\leqslant i<k.$ Then the pebbling configuration is $r$-solvable if a
vertex $r$ can be caught. For example, on Figure \ref{policepath} it is
possible to pebble from the vertex $c$ (cop) to the vertex $r$ (robber). 

Any
pebbling configuration $C$ determines the \emph{citizen subgraph}%
\textit{\ }$G_C$ of $K_{n}^{2}$ induced by the vertex set $%
\mathbf{P\cup T.}$ The edge set of $G_C$ is determined by the
vertices that see each other. Any component of $G_C$
containing two or more cops we call a \emph{police component}\textit{. }
\begin{claim}
\label{claimpolicecomp} Any configuration whose citizen subgraph has a police
component is solvable.
\end{claim}
\begin{proof}
Let us consider a police component with vertices $v_{1},\ldots
,v_{k} $ such that $v_{1},v_{k}\in \mathbf{P},$ $v_{1},\ldots ,v_{k-1}\in 
\mathbf{T} $ and $v_{i-1}$ sees $v_{i}$ for $1\leqslant i<k$. We now use the
following strategy.  Without loss of generality, we
assume that $v_{1}$ and $v_{2}$ are in the same row $v_{1}\Box K_{n}$. 
Then we find vertex $%
v_1^{\prime }$ which is the intersection of $K_{n}\Box r$ and $v_{1}\Box K_{n}$
and make $v_1^\prime$ a citizen by moving a pebble from $v_{1}$. 
Now $r=v_0$ can be caught by $v_1^\prime, v_2,\ldots,v_k$.
$\hfill$
\end{proof}
Another sufficient condition for a pebbling configuration to be solvable is
the existence of a ``robocop'', a vertex with $4$ or more pebbles on it. In
that case any robber $r$ can be caught by sending two pebbles to either $%
K_{n}\Box r$ or $r\Box K_{n}$, making a cop there and moving a pebble to $r$
from this new cop. In Section \ref{comps} we prove that the probability that
such a ``robocop'' exists tends to zero. 
Hence, our goal is to
prove that almost every configuration of asymptotically more than $n$
pebbles on $K_{n}^{2}$ has a police component.

The next argument transforms the original problem of the solvability of a
pebbling configuration on $K^{2}$ to connectedness properties of a related
bipartite multigraph $B_{n,n}^{\prime }$. First, we observe that $K_{n}^{2}$
is isomorphic to the line graph of the complete bipartite graph $K_{n,n}$.
Indeed, both vertex sets are isomorphic to $\{1,\ldots,n\}^2$, and both
edge sets are isomorphic to pairs from $\{1,\ldots,n\}^2$ that share a
coordinate.
Similarly we construct a bipartite graph $B_{n,n}$ whose line graph is
isomorphic to $G_C$. The bipartite multigraph $%
B_{n,n}^{\prime }$ is constructed from $B_{n,n}$ by adding multiple edges
according to the multiplicity of pebbles on the vertices of $G_C$.
In other words, for every vertex $(u,v)\in G_C$ we place the edge 
$uv\in B_{n,n}^{\prime }$ with multiplicity $C((u,v))$.

\section{\label{secmodeldesc}Model Descriptions}

In this section we describe three different models for random bipartite
graphs and multigraphs . We compare them and determine asymptotic
implications from one to another which we can apply then to the pebbling
threshold on the product of graphs. In particular, we will be interested in
the property of having a large component (which will be shown to be a
police component almost surely).

The first model (Model A) is an analogue of the probabilistic model for
random graphs. In model A edges between any two vertices in different parts
of $B_{n,n}$ are mutually independent and have the same probability $p$.
Computations are easiest in this model, in which all graphs are simple.
The probability space corresponding to model A we
denote by $\mathcal{B}_{n,p}$.

The second model (Model B) is an analogue of the
static model for random graphs. First, the set of all bipartite simple
graphs on $n$ by $n$ vertices is denoted $\mathcal{B}(n)$. 
The set of graphs in $\mathcal{B}(n)$ with $M$ edges we denote by $\mathcal{B}%
(n,M)$. \noindent This model consists of $|\mathcal{B}(n,M)|=\binom{N}{M}$
different graphs, where $N=n^{2}$. Clearly, $\mathcal{B}(n)=
\cup_{M=0}^N\mathcal{B}(n,M).$

Finally, we need a generalization of the second model for the case of
bipartite multigraphs (Model B$^{\prime }$). As it was defined in the
previous section, the edges of the line graph represent pebbles; therefore
we need a multiple edge model to reflect this situation. We denote by $%
\mathcal{B}^{\prime }(n,m)$ the set of all bipartite multigraphs on $n$ by $n$
vertices with $m$ edges$.$ \noindent Model B$^{\prime }$ consists of
precisely $|\mathcal{B}^{\prime }(n,m)|=
\QATOPD\langle \rangle {N}{m}$ different graphs, where $N=n^2$. 
Finally, we define $%
\mathcal{B}^{\prime }(n)=\cup_{m=0}^\infty \mathcal{B}%
^{\prime }(n,m).$

The multiple edge model for random graphs was considered in \cite{austin}.
It was shown that the differences between simple graphs and multigraphs are
negligible in most cases. Janson et al. \cite{birth} give a detailed
analysis of the multigraph model using an algebraic approach. 
We are going to show that, for the right translation of parameters, certain
properties that hold in model A will transfer to hold in B, and then to
B$^\prime$ as well.

\section{Connections Between Models}

Models A and B are very closely related to each other, provided that $M$ is
about $pN$, which is the expected number of edges of a graph in $%
\mathcal{B}_{n,p}$. In fact, these two models are asymptotically equivalent
to each other for any convex property. Call a family of multisets $\mathcal M$
\emph{increasing}\ if $A\subseteq B$\ and $A\in \mathcal{M}$ implies that $%
B\in \mathcal{M}$, \emph{decreasing}\ if $%
A\subseteq B$\ and $B\in \mathcal{M}$ implies that $A\in \mathcal{M}$. 
A family which is
either increasing or decreasing is called \emph{monotone}. Finally, a family 
$\mathcal{M}$ is \emph{convex}\ if $A\subseteq B\subseteq C$
and $A,C\in \mathcal{M}$ imply that $B\in \mathcal{M}$. 
Also, given a property $S$ we shall say that \emph{almost every} 
(\emph{a.e.}) graph in the probability space $\mathcal{M}$
has property $S$ if $\Pr [G\in \mathcal{M}:G$ has $\mathbf{S}]\rightarrow 1,$
as $n\rightarrow \infty .$

The equivalence of models A and B follows from the general equivalence of
the probabilistic and static models in random graphs, which was proven by 
Bollob\'{a}s (see \cite{bollobasgraph,palmer}).  Here we state the result
for random bipartite graphs.
\begin{result}
\label{equivalence} Let $N=n^2$ and let $0<p=p(n)<1$ be such that 
$pN\rightarrow \infty $
and $(1-p)N\rightarrow \infty $ as $n\rightarrow \infty ,$ and let $%
\mathbf{S}$ be a property of graphs.
\begin{enumerate}
\item  Suppose that $\varepsilon >0$ is fixed and that a.e. graph in $%
\mathcal{B}(n,M)$ has $\mathbf{S}$ whenever 
\begin{equation*}
(1-\varepsilon )pN<m<(1+\varepsilon )pN\ .
\end{equation*}
Then a.e. graph in $\mathcal{B}_{n,p}$ has $\mathbf{S}$.
\item  If $\mathbf{S}$ is a convex property and a.e. graph in $\mathcal{B}%
_{n,p}$ has $\mathbf{S}$, then a.e. graph in $\mathcal{B}(n,M)$ has $\mathbf{%
S}$ for $M=\left\lfloor pN\right\rfloor $.
\end{enumerate}
\end{result}
Bollob\'{a}s' technique is on the boolean lattice applied to graphs, so we
can apply it to bipartite graphs equally well since we are still considering
the boolean algebra in models A and B.
Next we establish a relationship between models B and B$^\prime$.

The \emph{support} of multigraph $G\in
\mathcal{B}^{\prime }(n)$\ is the simple graph obtained by identifying the
parallel edges of $G$.\ We denote the support by $\Lambda _{G}$. Obviously, $%
\Lambda _{G}\in \mathcal{B}(n)$. We call the number of edges in the support
the \emph{size} of the support of $G$, written $Z=Z_G= ||\Lambda _{G}||$. 
(Here we use the notation $||\cdot ||$ because we are counting edges rather
than vertices.)
The set of all graphs $G\in \mathcal{B}^{\prime }(n,m)$\ with the
same support size $Z_{G}=s$\ we denote $\Lambda (n,m,k)$.

An equivalent setting for the last definition is to consider $m$ unlabeled
balls placed in $N$ distinct urns. Then for $N=n^{2}$ the set $\Lambda _{G}$
represent the set of non-empty urns and $\Lambda (n,m,s)$ is the set of \
distributions into exactly $s$ of $N$ urns. We need to find the average size
of the support in this model. The probability that $G$ has support of size $%
s,$ for $0\leqslant s\leqslant m,$ is 
\begin{equation*}
\begin{tabular}{lll}
$\Pr [Z_{G}=s]$ & $=$ & $\dfrac{\binom{N}{s}\QATOPD\langle
\rangle {s}{m-s}}{\QATOPD\langle \rangle {N}{m}}$ \\ 
&  &  \\ 
& $=$ & $\dfrac{\binom{N}{s}\binom{m-1}{m-s}}{\binom{N+m-1}{m}}$ \\ 
&  &  \\ 
& $=$ & $\dfrac{\binom{N}{s}\binom{(N+m-1)-N}{m-s}}{\binom{N+m-1}{m}}$\ .
\end{tabular}
\end{equation*}
\newline
The last expression means that the random variable $Z=Z_{G}$
follows the hypergeometric distribution $\mathbf{H}$ with parameters $%
\mathbf{H}(N+m-1,N,m)$. The hypergeometric distribution $\mathbf{H}(L,k,l)$
describes the number of white balls in the sample of size $l$ chosen
randomly (without replacement) from an urn containing $L$ balls, of which $k$
are white and $L-k$ are black.
Direct computations give us the expected value and the variance of $%
Z_{G}$. Indeed, the general formula for $\mathbf{E}[Z^{k}]$ is 
\begin{equation}
\begin{tabular}{lll}
$\mathbf{E}[Z^{k}]$ & $=$ & $\overset{N}{\underset{s=0}{\sum }}s^{k}\Pr
[Z=s] $ \\ 
&  &  \\ 
& $=$ & $\overset{N}{\underset{s=0}{\sum }}
s^{k}\dfrac{\binom{N}{s}\binom{m-1}{m-s}}{\binom{N+m-1%
}{m}}$ \\ 
&  &  \\ 
& $=$ & $\dfrac{Nm}{N+m-1}\overset{N}{\underset{s=0}{\sum }}
s^{k-1}\dfrac{\binom{N-1}{s-1}\binom{%
m-1}{m-s}}{\binom{(N+m-1)-1}{m-1}}$ \\ 
&  &  \\ 
& $=$ & $\dfrac{Nm}{N+m-1}\overset{N-1}{\underset{r=0}{\sum }}
(r+1)^{k-1}\dfrac{\binom{N-1}{r}%
\binom{m-1}{(m-1)-r}}{\binom{N+m-2}{m-1}}$ \\ 
&  &  \\ 
& $=$ & $\dfrac{Nm}{N+m-1}E((Y+1)^{k-1})$\ ,%
\end{tabular}
\label{etok}
\end{equation}
where $Y$ is a hypergeometric random variable with parameters $\mathbf{H}%
(N+m-2,N-1,m-1)$. 

Setting $k=1$ in the last line of equation (\ref{etok}), we obtain 
\begin{equation}
\mathbf{E}[Z_{G}]=\dfrac{Nm}{N+m-1}=mq\ ,  
\label{expectedvalue}
\end{equation}
with $q=\dfrac{N}{N+m-1}$. The intuitive idea is that, for a large value of $%
N$, the average support size is close to $m.$ If the number of edges $%
m=m(n)\in o(N)$
then the value of $q$ is close to one. According to the second moment
method, if the variance of random variable $Z_{G}$ is relatively
small then the value of $Z_{G}$ almost always stays close to the
mean. Indeed, in equation (\ref{etok}) if $k=2$ then 
\begin{equation*}
\mathbf{E}[Z^{2}]=\dfrac{Nm}{N+m-1}\mathbf{E}[Y+1]=\dfrac{Nm}{N+m-1}\left( 
\dfrac{(N-1)(m-1)}{N+m-2}+1\right) .
\end{equation*}
Therefore, the variance is 
\begin{equation} 
\begin{tabular}{lll}
$\func{Var}[Z_{G}]$ & $=$ & $\mathbf{E}[Z^{2}]-(\mathbf{E}%
[Z])^{2} $ \\ 
&  &  \\ 
& = & $\dfrac{Nm}{N+m-1}\left( \dfrac{(N-1)(m-1)}{N+m-2}+1-\dfrac{Nm}{N+m-1}%
\right) $ \\ 
&  &  \\ 
& = & $mq\Big[ (N-1)(m-1)(N+m-1)+(N+m-2)(N+m-1)- $ \\
&  &  \\ 
& & \qquad $-Nm(N+m-2) \Big]\ \Big/\ (N+m-2)(N+m-1) $\\
&  &  \\ 
& = & $\dfrac{mq}{N+m-2}\left( \dfrac{Nm-N-m+1}{N+m-1}\right) $ \\
&  &  \\ 
& = & $\dfrac{mq}{N+m-2}(N-1)(1-q)\ .$ \\%
&  &  \\ 
\end{tabular}
\label{variance}
\end{equation}
For $m(n)\in o(N)$ we have from Equations (\ref{expectedvalue}) and (\ref
{variance}) that 
\begin{equation*}
\dfrac{\func{Var}[Z]}{(\mathbf{E}[Z])^{2}}=\dfrac{N-1}{N+m-2}\left( \dfrac{%
1-q}{mq}\right) \rightarrow 0\ ,
\end{equation*}
as $n\rightarrow \infty .$ Hence, $\func{Var}[Z_{G}]\in o((\mathbf{%
E}[Z_{G}])^{2}).$ 

Janson et al. \cite{randomjanson} suggested the following notation to
measure more precisely the closeness of a random variable to its mean.
\begin{notation}
Let $\{X_{n}\}_{n=1}^{\infty }$ be a sequence of random variables and $%
\{a_{n}\}_{n=1}^{\infty }$ a sequence of positive real numbers. We write 
\begin{equation*}
X_{n}=o_{p}(a_{n})
\end{equation*}
if, for every $\varepsilon >0,$ almost always $|X_{n}|<\varepsilon a_{n}$ $( 
$i.e. $\Pr [|X_{n}|<\varepsilon a_{n}]\rightarrow 1$ as $n\rightarrow \infty
).$
\end{notation}
This definition is analogous to $o(\cdot )$, but with probability involved.
\begin{lemma}\label{suppsize}
\label{suppism}Let $q=\dfrac{N}{N+m-1}$ and $m\in m(n)\subseteq o(N)$. Then 
\begin{equation*}
Z_{G}=mq+o_{p}(mq)\ ,
\end{equation*}
for a.e. graph $G$ in $\mathcal{B}^{\prime }(n,m)$.
\end{lemma}
\begin{proof}
We are going to use the second moment method with the random
variable $Z=Z_{G}$. We have $\mathbf{E}[Z]=mq$ by Equation \ref
{expectedvalue} and, using Equation \ref{variance} and Chebyshev's inequality, 
we obtain
\begin{equation*}
\Pr [|Z-mq|>\lambda ]\leqslant \dfrac{\sigma ^{2}}{\lambda ^{2}}=\dfrac{mq}{%
\lambda ^{2}}\dfrac{N-1}{N+m-2}(1-q)\leqslant \dfrac{mq}{\lambda ^{2}}\ 
\rightarrow 0\ 
\end{equation*}
for $\lambda =\varepsilon mq$ since $mq\rightarrow\infty$.
$\hfill$
\end{proof}
Now we are ready to establish the relationship between models B and B$%
^{\prime }$. The next theorem provides a criterion for any increasing
property that holds in $\mathcal{B}(n,M)$ to hold in $\mathcal{B}^{\prime
}(n,m)$ as well.
\begin{theorem}
\label{equiv} Let $\mathbf{S}$ be any increasing
property of graphs and $q=N/(N+m-1)$
for some $m\in o(N).$ Also, let $B_{\mathbf{S}}(n,M)\subseteq \mathcal{%
B}(n,M)$ and $B_{\mathbf{S}}^{\prime }(n,m)\subseteq \mathcal{B}^{\prime
}(n,m)$ denote those bipartite graphs and bipartite multigraphs,
respectively, having property $\mathbf{S}$. If for every sequence $M=M(n)$
such that $M=mq+o_{p}(mq)$ we have $\Pr [B_{\mathbf{S}}(n,M)]\rightarrow 1,$
as $n\rightarrow \infty ,$ then also $\Pr [B_{\mathbf{S}}^{\prime
}(n,m)]\rightarrow 1,$ as $n\rightarrow \infty $.
\end{theorem}
\begin{proof}
We are going to prove that $\Pr [B_{\mathbf{\bar{S}}}^{\prime
}(n,m)]\rightarrow 0.$ Let us consider the set 
\begin{equation*}
\mathbf{M}(\varepsilon )=\{M\text{ }|\text{ }|M-mq|\leqslant \varepsilon
mq\}\ ,
\end{equation*}
for some $\varepsilon >0.$ We assume in the hypothesis that for any $M\in 
\mathbf{M}(\varepsilon )$ we have
\begin{equation}
\Pr [B_{\mathbf{S}}(n,M)]\rightarrow 1\ ,  
\label{assum}
\end{equation}
whenever $n\rightarrow \infty .$ Then 
\begin{equation*}
\begin{tabular}{lll}
$\Pr [B_{\mathbf{\bar{S}}}^{\prime }(n,m)]$ & $=$ & $\sum\limits_{M\notin 
\mathbf{M}(\varepsilon )}\Pr [B_{\mathbf{\bar{S}}}^{\prime }(n,m)$ $|$ $%
Z_{G}=M]\Pr [Z_{G}=M]$ \\ 
&  &  \\ 
& $+$ & $\sum\limits_{M\in \mathbf{M}(\varepsilon )}\Pr [B_{\mathbf{\bar{S}}%
}^{\prime }(n,m)$ $|$ $Z_{G}=M]\Pr [Z_{G}=M]\ .$%
\end{tabular}
\end{equation*}
The first sum in the last expression can be bounded from above by 
\begin{equation*}
\sum\limits_{M\notin \mathbf{M}(\varepsilon )}\Pr [Z_{G}=M]\ =\ \Pr
[Z_{G}\notin M(\varepsilon )]\ ,
\end{equation*}
which tends to zero by Lemma \ref{suppism}. 
For every graph $B\in \mathcal{B}(n,M)$ there are 
$\QATOPD\langle \rangle {M}{m-M}$ multigraphs 
$B^\prime\in \mathcal{B}^{\prime }(n,m)$ with $%
\Lambda _{B^\prime}=B$.
Moreover, $\mathbf{S}$ is increasing.
Therefore we can give an upper bound for the second sum of 
\begin{equation}
\begin{tabular}{lll}
$\sum\limits_{M\in \mathbf{M}(\varepsilon )}\Pr [B_{\mathbf{\bar{S}}%
}^{\prime }(n,m)\ |\ Z_{B^\prime}=M]$ & $\leqslant$ & 
$\sum\limits_{M\in \mathbf{M}%
(\varepsilon )}\dfrac{\QATOPD\langle\rangle{M}{m-M}|B_{\mathbf{\bar{S}}}(n,M)|}
{\QATOPD\langle \rangle {N}{M}}$ \\ 
&  &  \\ 
& $=$ & $\sum\limits_{M\in \mathbf{M}(\varepsilon )}\dfrac{\QATOPD\langle
\rangle {M}{m-M}\binom{N}{m}}{\QATOPD\langle \rangle {N}{M}}\Pr [B_{\mathbf{%
\bar{S}}}(n,M)]$ \\ 
&  &  \\ 
& $\leqslant $ & $\Pr [B_{\mathbf{\bar{S}}}(n,M^{\ast })]\sum\limits_{M\in 
\mathbf{M}(\varepsilon )}\dfrac{\QATOPD\langle \rangle {M}{m-M}\binom{N}{m}}{%
\QATOPD\langle \rangle {N}{M}}\ ,$%
\end{tabular}
\label{secondterm}
\end{equation}
where $M^{\ast }$ is the element of $\mathbf{M}(\varepsilon )$ that
maximizes $\Pr [B_{\mathbf{\bar{S}}}(n,M)\mathbf{].}$ The sum in the last
expression is a partial sum of probabilities of a hypergeometric random
variable and, therefore, does not exceed 1. Hence, the last line in (\ref
{secondterm}) is bounded from above by $\Pr [B_{\mathbf{\bar{S}}}(n,M^{\ast
})]$, which goes to zero, as $n\rightarrow \infty $, by assumption (%
\ref{assum}). Thus, $\Pr [B_{\mathbf{\bar{S}}}^{\prime }(n,m)]\rightarrow 0$
as $n\rightarrow \infty $, and the statement of the theorem follows.
$\hfill$
\end{proof}
The particular increasing property in which we are most interested is that
of containing a large component, of size proportional to $2n$.
We will show that such a connected component is almost surely a police
component.

\section{\label{comps} Large Components and Police Components}

We first note that, almost surely, all cops have only two edges.
Recall that $B^\prime\in\mathcal{B}^\prime (n,m)$ is chosen uniformly 
at random, 
where $m=\varpi n$ and $\varpi\rightarrow\infty$ arbitrarily slowly as
$n\rightarrow\infty$.  
The probability that there exists a vertex with $k$ pebbles on it is at most
\begin{equation*}
\dfrac{n^{2}\QATOPD\langle \rangle {n^{2}}{\varpi n-k}}{%
\QATOPD\langle \rangle {n^{2}}{\varpi n}}\thicksim n^{2}\left( \dfrac{\varpi
n}{n^{2}+\varpi n}\right) ^{k}\thicksim n^{2}\left( \dfrac{\varpi }{n}%
\right) ^{k}.
\end{equation*}
For $k>2$ the last expression tends to zero as $n\rightarrow\infty$.

We use this fact to show that connected components of linear size have
many cops.
\begin{lemma}\label{copcomp}
Let $H$ be a connected component of size $\alpha 2n$ in 
$B^\prime\in\mathcal{B}^\prime (n,m)$, where $m=\varpi n$.
Then almost surely $H$ is a police component.
\end{lemma}
\begin{proof}
Let $x=x(n)$ be the excess of edges in $B^{\prime },$ namely $%
x(n)=||B^{\prime }||-||\Lambda _{B^{\prime }}||$. 
Since almost surely all cops have exactly two edges, the number of cops 
$s=s(n)$ in $B^\prime$ is almost always equal to the excess $x(n)$. 
Using Lemma \ref{suppsize} (with $q=N/(N+m-1)$) we compute
$$x(n) \ \sim \ m - mq \ = \ \varpi n (1-q) \ \sim \ \varpi^2 $$
almost surely.
Given that there are $\varpi^2$ cops in $B^\prime$, an upper bound of the probability that
$H$ has at most one cop is
\begin{eqnarray*}
\dfrac{ \dbinom{qn-\alpha 2n}{\varpi^2} 
+\alpha 2n\dbinom{qn-\alpha 2n}{\varpi^2-1} }{\dbinom{qn}{\varpi^2}}
& \lesssim & 
\left(\dfrac{qn-\alpha 2n}{qn}\right)^{\varpi^2} 
\left(1+\dfrac{2\alpha n\varpi^2}{qn-2\alpha n-\varpi^2}\right) \\
& & \\
& \lesssim & 
e^{-2\alpha\varpi^2/q}
\left(\dfrac{2\alpha\varpi^2}{q-2\alpha-\varpi^2/n}\right)\ .
\end{eqnarray*}
We may assume that $H$ is small, so if $4\alpha<q$ the last term is
at most $\varpi^2 e^{-\varpi^2/2}\rightarrow 0$ as $n\rightarrow\infty$.
$\hfill$
\end{proof}
Finally, we prove that there is a connected component of linear size in
$\mathcal{B}_{n,p}$.  The following theorem was proven in \cite{delavega}
for the random graph $\mathcal{G}_{n,p}$.  The proof involved analyzing
the hitting time of a certain parameter in a random walk and used no
special property of the graph structure.  Here we modify the result for
the random bipartite graph $\mathcal{B}_{n,p}$.  The same method yields
the following result, which we state without proof.
\begin{result}\label{dlv}
Let $\beta>\ln 16$, $p=\beta/n$, and $B_n\in\mathcal{B}_{n,p}$.
Then almost surely there is a path in $B_n$ of length at least
$(1-(\ln 16)/\beta)2n$.
\end{result}

\section{\label{secthrk2} Proof of Theorem \ref{thmk2}}

Now we prove that $\tau_{\mathcal{K}^2}=\Theta(\sqrt{N})$.
\begin{proof}
We recall that the pebbling threshold of every graph sequence is in 
$\Omega(\sqrt{N})$.  Therefore we need only show that
$\tau_{\mathcal{K}^2}=O(\sqrt{N})$.  Write $N=n^2$, let $m=\varpi n$, where
$\varpi=\varpi(n)\rightarrow\infty$ arbitrarily slowly, and let $C$ be a
randomly chosen configuration from $\mathbf{C}_{N,m}$.
Let $B^\prime_{n,n}\in \mathcal{B}^\prime (n,m)$ be the bipartite 
multigraph associated with $C$, and $B_{n,n}\in \mathcal{B}(n,M)$ 
be the simple bipartite graph determined by the support of $B^\prime_{n,n}$.
Lemma \ref{suppsize} implies $M=mq+o_p(mq)$, where $q=M/(M+n-1)$.

Let $p=M/N$ and consider the probability space $\mathcal{B}_{n,p}$.
For a graph $G$ let $\mathbf{S}=\mathcal{S}(G)$ be the property that $G$
has a connected component of size at least $\alpha |G|$, where
$\alpha=q/4$.  Define $\beta >\ln 16$ by $\alpha=1-(\ln 16)/\beta$ and let
$p^\prime=\beta/n$.  Then Result \ref{dlv} implies that almost every
graph in $\mathcal{B}_{n,p^\prime}$ has $\mathbf{S}$.  Since almost surely 
$p\sim mq/N\sim (\varpi/n)e^{-\varpi/n} > p^\prime$, and $\mathbf{S}$ is
an increasing property, almost every graph in $\mathcal{B}_{n,p}$ has
$\mathbf{S}$.
Every increasing property is also convex.  Thus Theorem
\ref{equivalence} assures that almost every graph in $\mathcal{B}(n,M)$
has $\mathbf{S}$.  Then Theorem \ref{equiv} implies that almost every
graph in $\mathcal{B}^\prime (n,m)$ has $\mathbf{S}$.  Let $H$ be such
a connected component of $B^\prime_{n,n}$ of size at least $\alpha 2n$.
According to Lemma \ref{copcomp} $H$ is almost surely a police component.
Finally, let $H_C$ be the corresponding connected component of the
citizen subgraph $G_C$ of the configuration $C$.  Since $H_C$ is a police
component, Claim \ref{claimpolicecomp} implies that $C$ is solvable.
This finishes the proof.
$\hfill$ 
\end{proof}

\section{Future Research}

Consider the graph $K_n^d=K_n^{d-1}\Box K_n$, and the sequence
$\mathcal{K}^d=\{K_1^d,\ldots,K_n^d,\ldots\}$.
If Conjecture \ref{hurlcon} is true then induction would show that
$\tau_{\mathcal{K}^d}=\Theta(\sqrt{N})$ for all $d$.  On the surface
such a result might be surprising, considering the sparcity of the graphs
(size $n^d$, degree $d(n-1)$).  However, its low diameter and high
structure make such a result believable.

Another interesting test for Conjecture \ref{hurlcon} is the sequence
of $n$-dimensional cubes $\mathcal{Q} =\{Q^1,\ldots,Q^n,\ldots\}$, where
$Q^n=Q^{n-1}\Box Q^1$, and $Q^1$ is the path on two vertices.  Because
$\mathcal{Q}^2$ is a subsequence of $\mathcal{Q}$, we must have
$\tau_{\mathcal{Q}^2}=\tau_\mathcal{Q}$.  Therefore, if $\tau_\mathcal{Q}
=\Theta(N^\alpha f(N))$ for some function $f(N)$, one can see that $f(N)$
must submultiplicative; i.e. $f(xy)\leqslant f(x)f(y)$ must hold.  
The best result 
to date is that $\tau_\mathcal{Q}\in \Omega(N^{1-\epsilon})\cap O(N/\lg N)$
for all $\epsilon>0$ (see \cite{cw}).

\section*{Acknowledgments}

The authors wish to thank Graham Brightwell for useful hallway conversations
during his visit to ASU.


\end{document}